\documentstyle{article}
\begin{document}
\renewcommand{\baselinestretch}{1.2}
\begin{center}
{\Large\bf {Linear mappings of local preserving-majorization\\ on matrix algebras}{\footnote{This work is supported by the
National Natural Science Foundation of China (No 10271012)
 and the Science Foundation of Hangzhou Dianzi
University}}} \vspace{0.4cm}
 \\
 {\bf Jun Zhu \footnote{E-mail address: zhu$\_$gjun@yahoo.com.cn}, Changping Xiong}\\
\vspace{0.3cm} \small{Institute of Mathematics, Hangzhou
Dianzi University, Hangzhou 310018, People's Republic of
China\\
}
\end{center}
\underline{~~~~~~~~~~~~~~~~~~~~~~~~~~~~~~~~~~~~~~~~~~~~~~~~~~~~~~~~~~~~~~~~~
~~~~~~~~~~~~~~~~~~~~~~~~~~~~~~~~~~~~~~~~~~~~~~~~~~~~~~~~~~~~~~~~~~~~~~~~~~~}
\vspace{0.3cm}\ \vspace{0.2cm} {\bf{Abstract}}^^L {Let
$\textbf{M}_{n\times n}$ be the algebra of all $n\times n$ matrices.
For $x,y\in {\textbf{R}}^{n}$ it is said that $x$ is majorized by $y$ if there is a double stochastic matrix $A\in {\textbf{M}}_{n\times n}$ such that $x=Ay$ (denoted by $x\prec y$). Let $x=(x_{1},x_{2},\cdots,x_{n})$ and $x=(y_{1},y_{2},\cdots, y_{n})$ in $\textbf{R}^{n}$, then $x\prec y$ if and only if $\sum_{i=1}^{k}x_{i}\leq \sum_{i=1}^{k}y_{i}$ and $\sum_{i=1}^{n}x_{i}= \sum_{i=1}^{n}y_{i}$. Suppose that $\Phi$ is a linear mapping from ${\textbf{R}}^{n}$ into ${\textbf{R}}^{n}$, which is said to be strictly isotone if $\Phi(x)\prec \Phi(y)$ whenever $x\prec y$.
We say that an element $\alpha\in
{\textbf{R}}^{n}$ is a strictly all-isotone point if every strictly isotone $\varphi$ at
$\alpha$ (i.e. $\Phi(\alpha)\prec\Phi(y)$ whenever $x\in
{\textbf{R}}^{n}$ with $\alpha\prec x$, and $\Phi(x)\prec\Phi(\alpha)$ whenever $x\in
{\textbf{R}}^{n}$ with $x\prec \alpha$) is a strictly isotone. In this paper we
show that every $\alpha=(\alpha_{1},\alpha_{2},\cdots,\alpha_{n})\in {\textbf{R}}^{n}$ with $\alpha_{1}>\alpha_{2}>\cdots>\alpha_{n}$ is a strictly all-isotone point. \vspace{0.2cm}
\\
\vspace{0.2cm} {\it{AMS Classification}}: 15A04, 15A21, 15A51
}\\
\vspace{0.2cm}
\underline{~~~~~~~~~~~~~~~~~~~~~~~~~~~~~~~~~~~~~~~~~~~~~~~~~~~~~~~~~~~~~~~~~
~~~~~~~~~~~~~~~~~~~~~~~~~~~~~~~~~~~~~~~~~~~~~~~~~~~~~~~~~~~~~~~~~~~~~~~~~~~}
\vspace{0.1cm}\

\section*{1. Introduction and preliminaries}
~

Let ${\textbf{R}}^{n}$ be real n-dimensional Euclidean spaces.
$L({\textbf{R}}^{n},{\textbf{R}}^{m})$ stands for the set of all linear mappings from ${\textbf{R}}^{n}$
into ${\textbf{R}}^{m}$, and abbreviate $L({\textbf{R}}^{n},{\textbf{R}}^{n})$ to $L({\textbf{R}}^{n})$. We denotes by ${\textbf{M}}_{n\times m}$ the set of all $n\times m$ matrices. The cardinal number of a set $A$ is denoted by $\mid A\mid$. $\textbf{N}$ is
the set of non-negative integers.

For $x,y\in {\textbf{R}}^{n}$ it is said that $x$ is majorized by $y$ if there is a double stochastic matrix $A\in {\textbf{M}}_{n\times n}$ such that $x=Ay$. Let us write $x\sim y$ if $x\prec y$ and $y\prec x$. If we write $x=(x_{1},x_{2},\cdots,x_{n})^{T}$ and $y=(y_{1},y_{2},\cdots,y_{n})^{T}$, then $x\prec y$ if and only if
$\sum_{i=1}^{k}x_{i}\leq \sum_{i=1}^{k}y_{i}$ ($i=1,2,\cdots,n-1$) and $\sum_{i=1}^{n}x_{i}= \sum_{i=1}^{n}y_{i}$ (see [9, Theorem 11.2]). A linear mapping $\Phi\in L({\textbf{R}}^{n},{\textbf{R}}^{m})$ is said to be strictly isotone if $\Phi(x)\prec \Phi(y)$ whenever $x\prec y$; $\Phi$ is said to be a strictly isotone at $\alpha$ if $\Phi(P\alpha)\prec\Phi(x)$ whenever $x\in
{\textbf{R}}^{n}$ with $\alpha\prec x$, and $\Phi(x)\prec\Phi(\alpha)$ whenever $x\in
{\textbf{R}}^{n}$ with $x\prec \alpha$. A vector $\alpha\in {\textbf{R}}^{n}$ is said to be strictly all-isotone point if every strictly isotone at $\alpha$ is a strictly isotone.

Some of our notations and symbols are explained as the following.

$\textbf{R}^{n}$: the set of all $n\times 1$ real column vectors.

$\textbf{R}_{n}$: the set of all $n\times 1$ real row vectors.

$\textbf{P}_{n}$: the set of all $n\times n$ permutation matrices.

$\textbf{S}_{n}$: the set of all $n\times n$ double stochastic matrices.

With the development of majorization problem, preserving majorization have attracted much attention of
mathematicians as an active subject of research in algebras. We describe some of the results related to ours. T. Anto [1] obtained the following interesting result.

{\bf{Theorem 1.1}} A linear mapping $\Phi: {\textbf{R}}^{n}\rightarrow {\textbf{R}}^{n}$ satisfies $\Phi(x)\prec \Phi(y)$ whenever $x\prec y$ if and onli if one of the following holds:

1) $\Phi(x)=(tr x)a$ for some $a\in {\textbf{R}}^{n}$;

2) $\Phi(x)=\alpha P(x)+\beta (tr x)e$ for some $\alpha, \beta\in {\textbf{R}}$ and $P\in \textbf{P}_{n}$.

Hereafter the linear preservers of majorization are fully characterized on the algebra of all $n\times n$ matrices by C.K. Li and E. Poon in [7].
A. Armandnejad, F. Akbarzadeh and Z. Mohammadi [2] obtained many results of preserving row and column-majorization on ${\textbf{M}}_{n\times m}$.
On the other hand, over the past few
years a considerable attention has been paid to the question of determining derivations (multiplicative mappings) through
one point derivable (multiplicative) mappings (see [3,5,6,8,10-12] and references therein). The purpose of this paper is to characterize those linear mappings of preserving-majorization at one fixing point.

This paper is organized as follows: in section 2, first we will introduce some notations of preserving majorization at one point, then we will give the main theorem in this paper and two lemmas which requires them in the proof of the main theorem. In section 3 we will give the proof of the main theorem.

\section*{2. Two lemmas and the main theorem}
~

Given two real vectors $x=(x_{1},x_{2},\cdots,x_{n})^{T}$ and $y=(y_{1},y_{2},\cdots,y_{n})^{T}$ in ${\textbf{R}}^{n}$, let $x^{\cdot}=(x_{1}^{\cdot},x_{2}^{\cdot},\cdots,x_{n}^{\cdot})^{T}$ and $x_{\cdot}=(x_{\cdot 1},x_{\cdot 2},\cdots,x_{\cdot n})^{T}$ denote $x$ and $y$ decreasing rearrangement and increasing rearrangement, respectively, i.e. $x_{1}^{\cdot}\geq x_{2}^{\cdot}\geq\cdots\geq x_{n}^{\cdot}$ and $x_{\cdot 1}\leq x_{\cdot 2}\leq\cdots\leq x_{\cdot n}$. The trace of $x$ is $tr(x)=\sum_{k=1}^{n}x_{k}$.
The set $\{P\in \textbf{P}_{n}: (P^{T} x)^{T}  y^{\cdot}=M(x,y)\}$ and $\{P\in \textbf{P}_{n}: (P^{T} x)^{T} y^{\cdot}=m(x,y)\}$ are denoted by $I_{M}(x,y)$ and $I_{m}(x,y)$, respectively, where $M(x,y)=(x^{\cdot})^{T}y^{\cdot}$ and $m(x,y)=(x_{\cdot})^{T}y^{\cdot}$.
We may identify an $n\times n$ matrix $A_{\Phi}$ as a linear mapping $\Phi$ on $\textbf{R}^{n}$.

{\bf{Definition 2.1}} {\it Let $\Phi: {\textbf{R}}^{n}\rightarrow {\textbf{R}}^{n}$ be a linear mapping.

1) $\Phi$ is said to be a strictly left-isotone at $\alpha\in {\textbf{R}}^{n}$ if $\Phi(y)\prec \Phi(P\alpha)$ whenever $P\in \textbf{P}_{n}$ and $y\in  {\textbf{R}}^{n}$ with $y\prec \alpha$.

2) $\Phi$ is said to be a strictly right-isotone at $\alpha\in {\textbf{R}}^{n}$ if $\Phi(P\alpha)\prec \Phi(y)$ whenever $P\in \textbf{P}_{n}$ and $y\in {\textbf{R}}^{n}$  with $\alpha\prec y$.

3) $\Phi$ is said to be a strictly isotone at $\alpha\in {\textbf{R}}^{n}$ if $\Phi(y)\prec \Phi(\alpha)$ whenever $y\in  {\textbf{R}}^{n}$ with $y\prec \alpha$, and $\Phi(\alpha)\prec \Phi(y)$ whenever $y\in  {\textbf{R}}^{n}$ with $\alpha\prec y$.

4) $\Phi$ is said to be a strictly preserving equivalence at $\alpha\in {\textbf{R}}^{n}$ if $\Phi(y)\sim \Phi(\alpha)$ whenever $y\in  {\textbf{R}}^{n}$ with $y\sim \alpha$.}

In the rest part of this paper, we always assume that $\alpha=(\alpha_{1},\alpha_{2},\cdots,\alpha_{n})^{T}$ with $\alpha_{1}>\alpha_{2}>\cdots>\alpha_{n}$. The following theorem is our main result.

{\bf{Theorem 2.2}} {\it Let $\Phi: {\textbf{R}}^{n}\rightarrow {\textbf{R}}^{n}$ be a linear mapping. Then the following statements are mutually equivalent:

1) $\Phi$ is a strictly left-isotone at $\alpha\in {\textbf{R}}^{n}$.

2) $\Phi$ is a strictly right-isotone at $\alpha\in {\textbf{R}}^{n}$.

3) $\Phi$ is a strictly isotone at $\alpha\in {\textbf{R}}^{n}$.

4) $\Phi$ is a strictly preserving equivalence at $\alpha\in {\textbf{R}}^{n}$.

5) $\Phi$ is a strictly isotone.}

The following most famous inequality for vectors is due to G.H. Hardy, J.E. Littlewood, G. Polya in [4, Theorem 368].

{\bf{Lemma 2.3}} {\it Let $x=(x_{1},x_{2},\cdots,x_{n})^{T}$ and $y=(y_{1},y_{2},\cdots,y_{n})^{T}$ be two vectors in ${\textbf{R}}^{n}$, then

(1) $m(x,y)=(x_{\cdot})^{T}y^{\cdot}\leq (P x)^{T}y\leq (x^{\cdot})^{T}y^{\cdot}=M(x,y), \forall P\in \textbf{P}_{n}.$

(2) If $y^{\cdot}=(y_{1},y_{2},\cdots, y_{n})$ with $y_{1}>y_{2}>\cdots>y_{n}$, then $(P x)^{T}y^{\cdot}=M(x,y)$ if and only if $Px=x^{\cdot}$.

(3) If $y^{\cdot}=(y_{1},y_{2},\cdots, y_{n})$ with $y_{1}>y_{2}>\cdots>y_{n}$, then $(P x)^{T}y^{\cdot}=m(x,y)$ if and only if $Px=x_{\cdot}$.}

{\bf Proof.} (1) The inequality for vectors is due to G.H. Hardy, J.E. Littlewood, G. Polya in [4, Theorem 368].

(2) The fact that the condition $Px=x^{\cdot}$ is sufficient by the inequality for vectors in (1). We only need to prove the necessity of the statement.
In fact, if $Px\neq x^{\cdot}$, we write $Px=(x_{P(1)},x_{P(2)},\cdots,x_{P(n)}),$ then there are $1\leq m<k\leq n$ such that
$x_{P(k)}<x_{P(m)}$. Since $(x_{P(m)}-x_{P(k)})(y_{k}-y_{m})>0$, we have $x_{P(k)}y_{k}+x_{P(m)}y_{m}<x_{P(m)}y_{k}+x_{P(k)}y_{m}$. It follows that
$$\begin{array}{rcl}(Px)^{T}y^{\cdot}&=&\sum_{i=1}^{n}x_{P(i)}y_{i}=\sum_{i\neq k,m}x_{P(i)}y_{i}+x_{P(k)}y_{k}+x_{P(m)}y_{m}\\&<&\sum_{i\neq k,m}x_{P(i)}y_{i}+x_{P(m)}y_{k}+x_{P(k)}y_{m}\leq M(x,y).\end{array}$$
Hence $(P x)^{T}y^{\cdot}\neq M(x,y)$.

(3) The statement can be proved by imitating the proof in (2).

This completes the proof of the lemma. $\Box$

{\bf{Lemma 2.4}} {\it Let $x=(x_{1},x_{2},\cdots,x_{n})^{T}$ and $y=(y_{1},y_{2},\cdots,y_{n})^{T}$ be two vectors in ${\textbf{R}}^{n}$ and $y_{1}>y_{2}>\cdots>y_{n}$. If there are $k$ unequal components at least in all entries of $x$,
then $\mid I_{M}(x,y^{\cdot})\mid\leq (n-k+1)!$ and $\mid I_{m}(x,y^{\cdot})\mid\leq (n-k+1)!$}.

{\bf Proof.} For $P\in  I_{M}(x,y^{\cdot})$, i.e. $(P^{T} x)^{T}  y^{\cdot}=M(x,y)=(x^{\cdot})^{T}y^{\cdot}$, it implies from Lemma 2.3 that
$$P^{T} x=(x_{P^{T}(1)}, x_{P^{T}(2)},\cdots, x_{P^{T}(n)})^{T}=x^{\cdot},$$
i.e. $x_{P^{T}(1)}\geq x_{P^{T}(2)}\geq\cdots\geq x_{P^{T}(n)}$.
Note that there are $k$ unequal components at least in all entries of $x$, then we have
$$\begin{array}{rcl}&&x_{P^{T}(1)}=x_{P^{T}(2)}=\cdots=x_{P^{T}(l_{1})}\\&>&
x_{P^{T}(l_{1}+1)}=x_{P^{T}(l_{1}+2)}=\cdots=x_{P^{T}(l_{2})}\\&>&
x_{P^{T}(l_{2}+1)}=x_{P^{T}(l_{2}+2)}=\cdots=x_{P^{T}(l_{3})}\\& &\cdots
\\&>&
x_{P^{T}(l_{k-2}+1)}=x_{P^{T}(l_{k-2}+2)}=\cdots=x_{P^{T}(l_{k-1})}
\\&>&x_{P^{T}(l_{k-1}+1)}\geq x_{P^{T}(l_{k-1}+2)}\geq\cdots\geq x_{P^{T}(n)}.\end{array}$$
It is easy to see that $\mid I_{M}(x,y)\mid\leq l_{1}!(l_{2}-l_{1})!\cdots(n-l_{k-1})!\leq (n-k+1)!$.

Similarly, we can prove that $\mid I_{m}(x,y^{\cdot})\mid\leq (n-k+1)!$. This completes the proof. $\Box$

{\bf{Note.}} Obviously $I_{M}(x,y^{\cdot})=\{P: (P^{T} x)^{T}y^{\cdot}=M(x,y)\}=\{P: x^{T}Py^{\cdot}=M(x,y)\}$ and $I_{m}(x,y^{\cdot})=\{P: (P^{T}x)^{T}y^{\cdot}=m(x,y)\}=\{P: x^{T}Py^{\cdot}=m(x,y)\}$.

\section*{3. The proof of the main theorem}
~

{\bf{The proof of Theorem 2.2}} From Definition 2.1 we can easily prove that "$1) \Rightarrow 4)$", "$2) \Rightarrow 4$" and "$3) \Rightarrow 4$".
It is easy to see that 5) implies 1)-3). So we only need to show that $4)$ implies $5)$.

Suppose that $\Phi$ is a linear mapping from $\textbf{R}^{n}$ into itself which satisfies strictly preserving equivalence at $\alpha\in {\textbf{R}}^{n}$, we write the $n\times n$ matrix $A_{\Phi}=(a_{ij})$, and abbreviate $A_{\Phi}$ to $A$. Since there is a real number $\lambda\in \textbf{R}$ such that $A+\lambda J=(b_{ij})$ with $b_{ij}>0$ (where every entries of the $n\times n$ matrix $J$ is 1). It suffices to prove the theorem in the case of the matrix $A=(a_{ij})$ with $a_{ij}>0$. We write $A=(a_{1}/a_{2}/\cdots/a_{n})$ ($A=(a^{1}|a^{2}|\cdots|a^{n})$), where $a_{j}$ ($a^{j}$) is the $j$th row (column) of $A$.

For any $P\in \textbf{P}_{n}$, it follows from $P\alpha\sim\alpha$ that $A(P\alpha)\sim A(\alpha)$.
Thus we may write
$$\begin{array}{lll}(A\alpha)^{\cdot}=(M_{1},M_{2},\cdots,M_{n})^{T},\end{array}$$
and$$\begin{array}{lll}&&AP\alpha=(a_{1}P\alpha,a_{2}P\alpha,\cdots,a_{n}P\alpha)^{T}\\
&=&(M_{k_{1}},M_{k_{2}},\cdots,M_{k_{n}})^{T}.\end{array}$$
Then there is an $a_{l}$ at least such that $a_{l}P\alpha=M_{1}=M(a_{l},\alpha)$, i,e, $P\in I_{M}(a_{l},\alpha)$.

We shall organize the proof of the theorem in a series of claims as follow.

{\bf Claim 1.} First, we show that $tr(a^{s})=tr(a^{t})$ for every $1\leq s,t\leq n$.

Taking two permutation matrices $P,Q\in \textbf{P}_{n}$, then $AP\alpha\sim A\alpha\sim AQ\alpha$ by $P\alpha \sim \alpha \sim Q\alpha$.
In particular $$\sum_{l=1}^{n}\sum_{k=1}^{n}a_{lk}\alpha_{P(k)}=\sum_{l=1}^{n}\sum_{k=1}^{n}a_{l}P\alpha=\sum_{k=1}^{n}M_{k}.$$
and $$\sum_{l=1}^{n}\sum_{k=1}^{n}a_{lk}\alpha_{Q(k)}=\sum_{l=1}^{n}\sum_{k=1}^{n}a_{l}Q\alpha=\sum_{k=1}^{n}M_{k}.$$
For any two nature numbers $1\leq s,t\leq n$ with $t\neq s$, we take a permutation matrix $P\in \textbf{P}_{n}$ such that $P(e_{s})=e_{t}$, $P(e_{t})=e_{s}$ and $P(e_{k})=e_{k}$ whenever $k\neq s,t$,
simultaneously we take $Q=I$ (where $I$ is the unit matrix in $\textbf{P}_{n}$).
Thus we have $\sum_{l=1}^{n}(a_{ls}\alpha_{t}+a_{lt}\alpha_{s})=\sum_{l=1}^{n}(a_{ls}\alpha_{s}+a_{lt}\alpha_{t})$, i.e. $\sum_{l=1}^{n}a_{ls}(\alpha_{t}-\alpha_{s})=\sum_{l=1}^{n}a_{lt}(\alpha_{t}-\alpha_{s})$.
Hence $tr(a^{s})=tr(a^{t})$.

{\bf Claim 2.} Suppose that every row $a_{i}$ of $A$ contains two unequal components at least in all entries of $a_{i}$.

We show that $A R=\left(
                  \begin{array}{cccc}
                    \gamma & \lambda & \cdots & \lambda \\
                    \lambda & \gamma & \cdots & \lambda \\
                    \cdots & \cdots & \cdots & \cdots \\
                    \lambda & \lambda & \cdots & \gamma \\
                  \end{array}
                \right)$
for some $R\in \textbf{P}_{n}$ and $\lambda,\gamma\in \textbf{R}$ with $\lambda\neq\gamma$.

Suppose that $n=2$. Then we may write
$A=\left(
     \begin{array}{cc}
       \lambda_{1} & \gamma_{1} \\
       \lambda_{2} & \gamma_{2} \\
     \end{array}
   \right)
$
and $\lambda_{1}+\lambda_{2}=\gamma_{1}+\gamma_{2}$ by Claim 1.
Without loss of generality, we assume that $\gamma_{1}>\lambda_{1}$.
Since $\left(
         \begin{array}{cc}
           0 & 1 \\
           1 & 0 \\
         \end{array}
       \right)\in \textbf{P}_{2}
$ and
$AP\alpha\sim A\alpha$, $\left(
                                                                 \begin{array}{cc}
                                                                    \gamma_{1} & \lambda_{1} \\
                                                                   \gamma_{2} & \lambda_{2} \\
                                                                 \end{array}
                                                               \right)\left(
                                                                        \begin{array}{c}
                                                                          \alpha_{1} \\
                                                                          \alpha_{2} \\
                                                                        \end{array}
                                                                      \right)
                                                               \sim\left(
     \begin{array}{cc}
       \lambda_{1} & \gamma_{1} \\
       \lambda_{2} & \gamma_{2} \\
     \end{array}
   \right)\left(
                                                                        \begin{array}{c}
                                                                          \alpha_{1} \\
                                                                          \alpha_{2} \\
                                                                        \end{array}
                                                                      \right)
$.
It follows that $\gamma_{1}\alpha_{1}+\lambda_{1}\alpha_{2}=\gamma_{2}\alpha_{2}+\lambda_{2}\alpha_{1}=M_{1}$ and $\gamma_{2}\alpha_{1}+\lambda_{2}\alpha_{2}=\gamma_{1}\alpha_{2}
+\lambda_{1}\alpha_{1}=M_{2}$.
Simple computation we obtain $\gamma_{2}=\lambda_{1}$, and $\gamma_{1}=\lambda_{2}$ by $\lambda_{1}+\lambda_{2}=\gamma_{1}+\gamma_{2}$.
Hence $A=\left(
     \begin{array}{cc}
       \gamma_{2} & \lambda_{2} \\
       \lambda_{2} & \gamma_{2} \\
     \end{array}
   \right).$

Suppose that $n\geq 3$. We divided the proof into four steps.

{\bf Step 1.}
We claim that every row $a_{i}$ of $A$ contains only two unequal components in all entries of $a_{i}$.

If the claim is not true, it follows that there is some row $a_{m}$ of $A$ such that it contains at least three unequal components in all entries of $a_{m}$.
By Lemmma 2.4, $\mid I_{M}(a_{m}, \alpha)\mid\leq (n-2)!$, and $\mid I_{M}(a_{k}, \alpha)\mid\leq (n-1)!$ whenever $k\neq m$.
Thus we have $$\sum_{k=1}^{n}\mid I_{M}(a_{k}, \alpha)\mid\leq (n-2)!+(n-1)(n-1)!<n!.$$
On the other hand, for every $P\in \textbf{P}_{n}$, we have $AP\alpha\sim A\alpha$. Thus there is some $a_{k}$ such that $a_{k}P\alpha=M(a_{k},\alpha)$, i.e. $P\in I_{M}(a_{k},\alpha)$.
It follows that  $$\sum_{k=1}^{n}\mid I_{M}(a_{k}, \alpha)\mid\geq\mid \textbf{P}_{n}\mid\geq n!,$$
which is a contradiction. Hence the claim holds.

{\bf Step 2.}
We claim that every components of $a_{i}$ is equal to same real number only except one, $\sum_{k=1}^{n}\mid I_{M}(a_{l},\alpha)\mid=n!$ and $I_{M}(a_{l},\alpha)\cap I_{M}(a_{k},\alpha)=\emptyset$ for $l\neq k$.

If the claim is not true, then there are an $a_{l_{0}}$, $P\in \textbf{P}_{n}$ and $2\leq k\leq n-2$ such that $a_{l_{0}}=(a,\cdots,a,b,\cdots,b)P$, where $a$ and $b$ appears $k$ times and $n-k$ times in components of $a_{l_{0}}$, respectively.
By Lemma 2.3, we have $\mid I_{M}(a_{l_{0}},\alpha)\mid\leq k!(n-k)!<(n-1)!$.
Note that $\mid I_{M}(a_{l},\alpha)\mid\leq (n-1)!$ whenever $l\neq l_{0}$, thus we obtain
$\sum_{l=1}^{n}\mid I_{M}(a_{l},\alpha)\mid<n!$.
On the other hand, for every $P\in \textbf{P}_{n}$, we have $AP\alpha\sim A\alpha$. Thus there is some $a_{k}$ such that $a_{k}P\alpha=M(a_{k},\alpha)$, i.e. $P\in I_{M}(a_{k},\alpha)$.
It follows that  $$\sum_{k=1}^{n}\mid I_{M}(a_{k}, \alpha)\mid\geq \mid \textbf{P}_{n}\mid\geq n!,$$ which is a contradiction.
Hence $$a_{l}=(\lambda_{l},\cdots,\lambda_{l},\gamma_{l},\lambda_{l},\cdots,\lambda_{l})$$ with $\lambda_{l}\neq\gamma_{l}$, $(l=1,2,\cdots,n).$

It is easy to see from Lemma 2.3 and the form of the above vector that $\mid I_{M}(a_{l},\alpha)\mid=(n-1)!$. So $\sum_{k=1}^{n}\mid I_{M}(a_{l},\alpha)\mid=n!$.
Since $\mid \textbf{P}_{n}\mid=n!$, we have $I_{M}(a_{l},\alpha)\cap I_{M}(a_{k},\alpha)=\emptyset$ for $l\neq k$.

{\bf Step 3.} We claim that every column $a^{s}$ of $A$ contains two $\gamma_{t}$'s at most in its components.

If the claim is not true, then there are three rows $a_{l}, a_{k}$ and $a_{m}$ of $A$ such that
$$\begin{array}{lll}a_{l}&=&(\lambda_{l},\cdots,\lambda_{l},\gamma_{l},\lambda_{l},\cdots,\lambda_{l}),
\\a_{k}&=&(\lambda_{k},\cdots,\lambda_{k},\gamma_{k},\lambda_{k},\cdots,\lambda_{k}),
\\a_{m}&=&(\lambda_{m},\cdots,\lambda_{m},\gamma_{m},\lambda_{m},\cdots,\lambda_{m}),\end{array}$$
where the $j$th components of $a_{l}$, $a_{k}$ and $a_{m}$ are $\gamma_{l}$, $\gamma_{k}$ and $\gamma_{m}$, respectively.
There are several possibilities.

{\it{Case 1. }} Suppose that $\gamma_{l}>\lambda_{l}$ and $\gamma_{k}>\lambda_{k}$.
There are $P\in \textbf{P}_{n}$ such that $a_{l}P\alpha=M(a_{l},\alpha)$, i.e. $P\in I_{M}((a_{l},\alpha)$.
Obviously $P\in I_{M}(a_{k},\alpha)$, this is a contradiction with $I_{M}(a_{l},\alpha)\cap I_{M}(a_{k},\alpha)=\emptyset$.

{\it{Case 2. }} Suppose that $\gamma_{l}<\lambda_{l}$ and $\gamma_{k}<\lambda_{k}$.
There are $P\in \textbf{P}_{n}$ such that $a_{l}P\alpha=M(a_{l},\alpha)$, i.e. $P\in I_{M}((a_{l},\alpha)$.
Obviously $P\in I_{M}(a_{k},\alpha)$, this is a contradiction with $I_{M}(a_{l},\alpha)\cap I_{M}(a_{k},\alpha)=\emptyset$.

{\it{Case 3. }} Suppose that $\gamma_{m}<\lambda_{m}$ and $\gamma_{k}<\lambda_{k}$.
There are $P\in \textbf{P}_{n}$ such that $a_{m}P\alpha=M(a_{m},\alpha)$, i.e. $P\in I_{M}((a_{m},\alpha)$.
Obviously $P\in I_{M}(a_{k},\alpha)$, this is a contradiction with $I_{M}(a_{l},\alpha)\cap I_{M}(a_{k},\alpha)=\emptyset$.

The other cases can beget the same contradiction by imitating the above Case 1, Case 2 and Case 3.
Thus every column $a^{k}$ of $A$ contains two $\gamma_{i}$'s at most in its components.

{\bf Step 4.}  We claim that every column $a^{s}$ of $A$ includes only one $\gamma_{t}$'s in its components.

If the claim is not true, then there are two rows $a_{k}$ and $a_{m}$ of $A$ such that
$$\begin{array}{lll}a_{k}&=&(\lambda_{k},\cdots,\lambda_{k},\gamma_{k},\lambda_{k},\cdots,\lambda_{k}),
\\a_{m}&=&(\lambda_{m},\cdots,\lambda_{m},\gamma_{m},\lambda_{m},\cdots,\lambda_{m}),\end{array}$$
where the $j$th components of $a_{k}$ and $a_{m}$ are $\gamma_{k}$ and $\gamma_{m}$, respectively.
At this time there must be some column $a^{l}$ of $A$ such that $a^{l}=(\lambda_{1},\lambda_{2},\cdots,\lambda_{n})^{T}$ and $tr(a^{l})=\sum_{i=1}^{n}\lambda_{i}$.

There are several possibilities.

{\it{Case 1. }} Suppose that there is some column $a^{j}$ of $A$ such that $a^{j}$ includes only one $\gamma_{t}$'s in its components.
By Claim 1, we have $\gamma_{t}+\sum_{i\neq t}\lambda_{i}=tr(a^{j})=tr(a^{l})=\sum_{i=1}^{n}\lambda_{i}$, i.e.
$\gamma_{t}=\lambda_{t}$, this is a contradiction with $\gamma_{t}\neq\lambda_{t}$.

{\it{Case 2. }} Suppose that every column $a^{j}$ of $A$ includes two $\gamma_{t}$'s or not in its components.
In this case, $n$ must be an even, which implies $n\geq 4$.
Without loss of generality, we write
$$A=\left(
\begin{array}{ccccc}
                    \gamma_{1} & \lambda_{1} & \lambda_{1}&\cdots & \lambda_{1} \\
                    \gamma_{2} & \lambda_{2} & \lambda_{2}&\cdots & \lambda_{2} \\
                    \lambda_{3}& \gamma_{3} & \lambda_{3}&\cdots & \lambda_{3} \\
                    \lambda_{4}& \gamma_{4} & \lambda_{4}&\cdots & \lambda_{4} \\
                    \cdots & \cdots & \cdots & \cdots & \cdots\\
                    \lambda_{n} & \lambda_{n} &\cdots& \cdots & \cdots \\
                  \end{array}
                \right).$$
Since $tr(a^{1})=tr(a^{2})=\sum_{i=1}^{n}\lambda_{i}$, $\gamma_{1}+\gamma_{2}=\lambda_{1}+\lambda_{2}$ and $\gamma_{3}+\gamma_{4}=\lambda_{3}+\lambda_{4}$.
Thus there are $0\neq c, d\in \textbf{R}$
such that $\gamma_{1}=\lambda_{1}+c$, $\gamma_{2}=\lambda_{2}-c$, $\gamma_{3}=\lambda_{3}+d$ and $\gamma_{4}=\lambda_{4}-d$.

Take $P_{0}=\left(
                                   \begin{array}{ccccc}
                                     0 & 0 & \cdots &0 &1\\
                                     1 & 0 & \cdots &0&0\\
                                     0 & 1 &\cdots &0&0\\
                                     \cdots&\cdots&\cdots&\cdots&\\
                                     0 & 0 & \cdots & 1&0\\
                                   \end{array}
                                 \right)\in \textbf{P}_{n}.
$ Then

a) Suppose that $c>0$ and $d>0$. Then $P_{0}\in I_{M}(a_{2},\alpha)\cap I_{M}(a_{3},\alpha)$ by Lemma 2.3, which is a contradiction with $I_{M}(a_{l},\alpha)\cap I_{M}(a_{k},\alpha)=\emptyset$ for $l\neq k$.

b) Suppose that $c>0$ and $d<0$. Then $P_{0}\in I_{M}(a_{2},\alpha)\cap I_{M}(a_{4},\alpha)$ by Lemma 2.3, which is a contradiction with $I_{M}(a_{l},\alpha)\cap I_{M}(a_{k},\alpha)=\emptyset$ for $l\neq k$.

c) Suppose that $c<0$ and $d>0$. Then $P_{0}\in I_{M}(a_{1},\alpha)\cap I_{M}(a_{3},\alpha)$ by Lemma 2.3, which is a contradiction with $I_{M}(a_{l},\alpha)\cap I_{M}(a_{k},\alpha)=\emptyset$ for $l\neq k$.

d) Suppose that $c<0$ and $d<0$. Then $P_{0}\in I_{M}(a_{1},\alpha)\cap I_{M}(a_{4},\alpha)$ by Lemma 2.3, which is a contradiction with $I_{M}(a_{l},\alpha)\cap I_{M}(a_{k},\alpha)=\emptyset$ for $l\neq k$.

Hence every column $a^{s}$ of $A$ includes only one $\gamma_{t}$'s in its components ($s=1,2,\cdots,n$).
Since $tr(a^{i})=tr(a^{j})(i,j=1,2,\cdots,n)$, we may write
$$A P=\left(
      \begin{array}{cccc}
        \lambda_{1}+\beta & \lambda_{1} & \cdots & \lambda_{1} \\
        \lambda_{2} & \lambda_{2}+\beta & \cdots & \lambda_{2} \\
         \cdots &  \cdots &  \cdots &  \cdots \\
        \lambda_{n} & \lambda_{n} & \cdots & \lambda_{n}+\beta \\
      \end{array}
    \right)
$$
for some $P\in\textbf{P}_{n}$ and $\beta\in \textbf{R}$.
Without loss of generality, we assume $\beta>0$. Since
$$(M_{1},M_{2},\cdots,M_{n})=(A\alpha)^{\cdot}\sim AP\alpha\sim AQ\alpha, \forall Q\in \textbf{P}_{n},$$
we have $M_{1}=(\lambda_{m}+\beta)\alpha_{1}+\sum_{k= 2}^{n}\lambda_{m}\alpha_{k}, (m=1,2,\cdots,n)$.
So $(\lambda_{m}+\beta)\alpha_{1}+\sum_{k= 2}^{n}\lambda_{m}\alpha_{k}=(\lambda_{l}+\beta)\alpha_{1}+\sum_{k= 2}^{n}\lambda_{l}\alpha_{k}$, i.e.
$\lambda_{m}=\lambda_{k}=\lambda$ for every $m,l=1,2,\cdots,n$, i.e.
 $$A P=\left(
                  \begin{array}{cccc}
                    \lambda+\beta & \lambda & \cdots & \lambda \\
                    \lambda & \lambda+\beta & \cdots & \lambda \\
                    \cdots & \cdots & \cdots & \cdots \\
                    \lambda & \lambda & \cdots & \lambda+\beta \\
                  \end{array}
                \right)$$
It is easy to verify that $A$ is a strictly isotone.

{\bf Claim 3.}
Suppose that there is some row $a_{m}$ of $A$ such that every components of $a_{m}$ is equal to same real number, i.e. $a_{m}=(\lambda_{m},\lambda_{m},\cdots,\lambda_{m})$. We claim that $a_{i}=(\lambda_{i},\lambda_{i},\cdots,\lambda_{i})$ for any $1\leq i\leq n$.

If the claim is not true, then there is some row $a_{k}$ of $A$ which contains two unequal components at least in all entries of $a_{k}$.
Without loss of generality we may assume that $a_{i}=(\lambda_{i},\lambda_{i},\cdots,\lambda_{i})$ whenever $1\leq i\leq k$, and $a_{i}$ contains two unequal components at least whenever $k+1\leq i\leq n$, i.e.
$$A=\left(
     \begin{array}{c}
       a_{1} \\
       a_{2} \\
       \cdots \\
       a_{k} \\
       A_{1} \\
     \end{array}
   \right),~~~~~ A_{1}=\left(
                    \begin{array}{c}
                      a_{k+1} \\
                      a_{k+2} \\
                      \cdots \\
                      a_{n} \\
                    \end{array}
                  \right),
$$
where every row of $A_{1}$ contains two unequal components at least. It is easy to verify from $A\alpha\sim AP\alpha\sim AQ\alpha$ that $A_{1}\alpha\sim A_{1}P\alpha\sim A_{1}Q\alpha$ for every $P,Q\in \textbf{P}_{n}$.
We write $(A_{1}\alpha)^{\cdot}=(N_{1},N_{2},\cdots,N_{n-k})^{T}$, i.e. $N_{1}\geq N_{2}\geq\cdots\geq N_{n-k}$. For every $P\in \textbf{P}_{n}$,
there is some row $a_{l} (k+1\leq l\leq n)$ of $A_{1}$ such that $a_{l}P\alpha=N_{1}$, i.e. $P\in I_{M}(a_{l},\alpha)$. Hence $\sum_{i=1}^{n-k}\mid I_{M}(a_{k+i},\alpha)\mid\geq n!$. On the other hand, $\mid I_{M}(a_{k+i},\alpha)\mid\leq (n-1)!$ by Lemma 2.4. Thus $\sum_{i=1}^{n-k}\mid I_{M}(a_{k+i},\alpha)\mid\leq (n-k)(n-1)!<n!$. This is a contradiction. Hence $a_{i}=(\lambda_{i},\lambda_{i},\cdots,\lambda_{i})$ for any $1\leq i\leq n$ or $$A=\left(
            \begin{array}{cccc}
              \lambda_{1} & \lambda_{1} & \cdots & \lambda_{1} \\
              \lambda_{2} & \lambda_{2} & \cdots & \lambda_{2} \\
              \dots &\dots &\dots &\\
             \lambda_{n} & \lambda_{n} & \cdots & \lambda_{n} \\
            \end{array}
          \right).
$$
It is easy to verify that $A$ is a strictly isotone. This completes the proof of the theorem.
$\Box$\\

{\bf{Corollary 3.1}} {\it Every $\alpha=(\alpha_{1},\alpha_{2},\cdots,\alpha_{n})\in {\textbf{R}}^{n}$ with $\alpha_{1}>\alpha_{2}>\cdots>\alpha_{n}$ is a strictly all-isotone point.}
\\~\\
\section*{Reference}\small
\begin{description}
  \item[1] T. Ando, Majorization, doubly stochastic matrices, and comparison of eigenvalues, Linear Algebra and its applications, 118 (1989) 163-248.
  \item[2] A. Armandnejad, F. Akbarzadeh, Z. Mohammadi, Row and column-majorization on $\textbf{M}_{n,m}$, Linear Algebra and its applications, 437 (2012) 1025-1032.
  \item[3] M.A. Chebotar, W.F. Ke, P.H. Lee, N.C.Wong, Mappings preserving zero products, Stud. Math. 155 (1) (2003) 77-94.
  \item[4] G.H. Hardy, J.E. Littlewood, G. Polya, Inequalites, Cambridge University Press, Cambridge, 1934.
  \item[5] J.C. Hou, X.F. Qi, Additive maps derivable at some points on J-subspace lattice algebras, Linear Algebra Appl. 429 (2008)
  1851-1863.
  \item[6] W. Jing, On Jordan all-derivable points of $B(H)$, Linear Algebra and Appl. 430(4) (2009) 941-946.
  \item[7] C.K. Li, Linear operators preserving directional majorization, Linear Algebra and its applications, 325 (2001) 141-146.
  \item[8] X.F. Qi, J.C. Hou, Linear maps Lie derivable at zero on J-subspace lattice algebras,
  Studia. Math. 197 (2010) 157-169.
  \item[9] J. Watrous, Theory of Quantum information, University of Waterloo, Waterloo, 2008.
  \item[10] J. Zhu, All-derivable points of operator algebras, Linear Algebra Appl. 427 (2007) 1-5.
  \item[11] J. Zhu, C.P. Xiong, All-derivable points in continuous nest algebras, J. Math. Anal. Appl. 340 (2008) 845-853.
  \item[12] J. Zhu, C.P. Xiong, H. Zhu, Multiplicative mappings at some points on matrix algebras, Linear Algebra Appl. 433 (2010) 914-927.

  \end{description}

\end{document}